\RequirePackage{fix-cm}
\documentclass[smallextended]{svjour3}       % onecolumn (second format)
\smartqed  % flush right qed marks, e.g. at end of proof
\usepackage{graphicx}
\usepackage{mathrsfs}
\usepackage{amssymb}
\usepackage{amsfonts}
\usepackage{amsmath}
\usepackage{epsfig}
\usepackage{color}
\usepackage{subfig}
\usepackage{epstopdf}
\usepackage{algorithm}
\usepackage{algpseudocode} 
\makeatletter
\def\BState{\State\hskip-\ALG@thistlm}
\makeatother
\newtheorem{assumption}{Assumption}
\begin{document}
	
	\title{Optimal On-Off Control for a Class of Discrete Event Systems
		with Real-Time Constraints\thanks{The authors' work is supported in part by a start-up funding provided by Middle Tennessee State University.}
	}
	
	
	\author{Lei Miao \and
			 Lijian Xu \and
			 Dingde Jiang         %etc.
	}
	
	
	\institute{Lei Miao \at
		Mechatronics Engineering \\
		1301 East Main Street, Box 19 \\ Murfreesboro, TN 37132-0001, USA \\
		Tel.: +1-615-898-2256\\
		Fax: +1-615-898-5697\\
		\email{lei.miao@mtsu.edu}             %\\		
	\and
		Lijian Xu\at
		Department of Electrical and Computer Engineering Technology \\
		Farmingdale State College, Farmingdale, NY 11735, USA \\
		\email{xul@farmingdale.edu}             %\\
	\and
		Dingde Jiang \at
		School of Astronautics and Aeronautic \\
		University of Electronic Science and Technology of China, Chengdu 611731, China\\
		\email{jiangdd@uestc.edu.cn}           %  \\
}
	
	\date{Received: date / Accepted: date}

	\maketitle

%%%%%%%%%%%%%%%%%%%%%%%%%%%%%%%%%%%%%%%%%%%%%%%%%%%%%%%%%%%%%%%%%%%%%%%%%%%%%%%%
\begin{abstract}

This paper studies an optimal ON-OFF\ control problem for a class of
discrete event systems with real-time constraints. Our goal is to minimize
the overall costs, including the operating cost and the wake-up cost,
while still guaranteeing the deadline of each individual task. In
particular, we consider the homogeneous case in which it takes the same
amount of time to serve each task and each task needs to be served by $d$
seconds upon arrival. The problem involves two subproblems:\ \textit{(i)} finding the
best time to wake up the system and \textit{(ii)} finding the best time to let the
system go to sleep. We study the two subproblems in both off-line and on-line settings. In the off-line case that all task information is known a priori, we combine sample path analysis and dynamic programming to come up with the optimal solution. In the on-line scenario where future task information is completely unknown, we show that the optimal time to wake up the system can be obtained without relying on future task arrivals. We also perform competitive analysis for on-line control and derive the competitive ratios for both deterministic and random controllers. 
\end{abstract}

%%%%%%%%%%%%%%%%%%%%%%%%%%%%%%%%%%%%%%%%%%%%%%%%%%%%%%%%%%%%%%%%%%%%%%%%%%%%%%%%
\section{Introduction}

There exists a large amount of Discrete Event Systems (DESs) that involve
allocation of resources to satisfy real-time constraints. One commonality of these DESs is that certain tasks must be completed by their deadlines in order to guarantee Quality-of-Service (QoS). Examples arise in wireless networks and computing systems, where communication and computing tasks must be transmitted/processed before the information they contain becomes obsolete \cite{MiaoMaoCGCTONS} \cite{Liu00}, and in manufacturing systems, where manufacturing tasks must be completed before the specified time in the production schedule \cite{PepCas00}. Another commonality of these DESs is that they all require the minimization of cost (e.g., energy). An interesting question then arises naturally:\ \emph{how can we allocate resources to such DESs so that the cost is minimized and the real-time constraints are also satisfied?} To answer this question, one often has to study the trade-off between minimizing the cost and satisfying the real-time constraints:\ processing the tasks at a higher speed makes it easier to satisfy the real-time constraints but harder to reduce the cost; conversely, processing the tasks at a lower speed makes it harder to satisfy the real-time constraints but easier to reduce the cost. This trade-off is often referred to as the energy-latency trade-off and has been widely studied in the literature \cite{MiaoMaoCGCTONS} \cite{GamNaPraUyZaInf02} \cite{ZaferTON09}.

In this paper, our objective is to utilize the energy-latency trade-off to minimize the cost while guaranteeing the real-time constraint for each task. Different from most existing papers that assume the system's service rate (the control variable)\ is a continuous function of time, we assume that the DES only operates at one of the two states: ON and OFF. One motivating example of such DES is wireless sensor networks, in which operation simplicity must be maintained. For example, the radio of a ZigBee wireless device can be either completely off or transmitting at a fixed-rate, e.g., 250kb/s in the 2.4GHz band. Another difference between this paper and others is that we assume that a wake-up cost is incurred whenever the system transits from the OFF state to the ON state.

In this paper, we solve both off-line and on-line optimal ON-OFF control problems. Our main contributions are two-fold: \textit{(i)} We combine \emph{sample path analysis} and \emph{Dynamic Programming} (DP) to obtain the optimal off-line solution and \textit{(ii)} We perform competitive analysis and derive the competitive ratios of both deterministic and random on-line controllers. Some results of this paper are previously shown in two conference papers: \cite{MiaoXuWTS2015} and \cite{Miao2017ACC_ON_OFF_Control}, which primarily focus on off-line control. One new contribution of this paper is the competitive analysis for on-line control. Another new contribution is that we introduce an idling cost in the system model. We point out that the addition of this idling cost makes our problem formulation more realistic because it often exists in real-world applications. For example, energy is consumed when a motor is spinning without any load attached and when a sensor is turned on, but not actively processing information. 

In this journal version, we improve some proofs to incorporate the idling cost; we also move all the proofs and tables to the appendix in order to enhance the continuity of the
analysis in the paper. The organization of the rest of the paper is as follows: in Section \ref{Sec_related_work}, we discuss related work; we introduce the system model and formulate our optimization problem in Section \ref{Sec_model}; the off-line and on-line results are presented in Sections \ref{Sec_offline} and \ref{Sec_online}, respectively; finally, we conclude in Section \ref{Sec_conclusions}.

\section{Related Work}
\label{Sec_related_work}
There are two lines of work that are closely related to this paper. One is transmission scheduling for wireless networks, in which the transmission rate of a wireless device is adjusted so as to minimize the transmission cost and satisfy real-time constraints. This line of work is initially studied in \cite{TON-Energy-Efficient_Wireless} with follow-up work in \cite{GamNaPraUyZaInf02} where a homogeneous case is considered, assuming all packets have the same deadline and number of bits. By identifying some properties of this convex optimization problem, Gamal et al. propose the ``MoveRight" algorithm in \cite{GamNaPraUyZaInf02} to solve it iteratively. However, the rate of convergence of the MoveRight algorithm is only obtainable for a special case of the problem when all packets have identical energy functions; in general the MoveRight algorithm may converge slowly. Zafer et al. \cite{Zafer09TAC} study an optimal rate control problem over a time-varying wireless channel, in which the channel state was modeled as
a Markov process. In particular, they consider the scenario that $B$ units of data must be transmitted by a common deadline $T,$ and they
obtain an optimal rate-control policy that minimizes the total energy
expenditure subject to short-term average power constraints. In \cite%
{Zafer07ITA} and \cite{Zafer08TIT}, the case of identical arrival time and individual deadline is studied by Zafer et. al. In \cite{NeedlyInfocom07}, the case of identical packet size and identical delay constraint is studied by Neely et. al. They extend the result for the case of individual packet size and identical delay constraint in \cite{NeedlyWirelessNetworks09}. In \cite{ZaferTON09}, Zafer et. al. use a graphical approach to analyze the case that each packet has its own arrival time and deadline. However, there are certain restrictions in their setting; for example, the packet that arrives later must have later deadlines. Wang and Li \cite{WangToWC2013} analyze scheduling problems for bursty packets with strict deadlines over a single time-varying wireless channel. Assuming slotted transmission and
changeable packet transmission order, they are able to exploit structural
properties of the problem to come up with an algorithm that solves the
off-line problem. In \cite{PoulakisTOVT2013}, Poulakis et. al. also study energy efficient scheduling problems for a single time-varying wireless channel. They consider a finite-horizon problem where each packet must be transmitted before $D_{\max }.$ Optimal stopping theory is used to find the optimal start transmission time between $[0,$ $D_{\max }]$ so as to minimize the expected energy consumption and the average energy consumption per unit of time. Zhong and Xu \cite{ZhongXuInfocom2008} formulated optimization problems that minimize the energy consumption of a set of tasks with task-dependent energy functions and packet lengths. In their problem formulation, the energy functions include both transmission energy and circuit power consumption. To obtain the optimal solution for the off-line
case with backlogged tasks only, they develop an iterative algorithm RADB whose complexity is $O(n^{2})$ ($n$ is the number of tasks). The authors show via simulation that the RADB algorithm achieves good performance when used in on-line scheduling. \cite{MiaoMaoCGCTONS} studies a transmission control problem for task-dependent cost functions and arbitrary task arrival time, deadline, and number of bits. They propose a GCTDA\ algorithm that solves the off-line problem efficiently by identifying certain critical tasks. The GCTDA algorithm is an extension to the CTDA algorithm \cite{MaoCTDA} designed by Mao and Cassandras for dynamic voltage scaling related applications. They extend the CTDA algorithm to multilayer scenarios in \cite{mao2014optimal}. Our model is different from all the above works by letting the system operate in one of the discrete modes and also including a wake-up cost at each time instant that the system transitions from OFF to ON state.

The other line of research studies On-OFF scheduling in Wireless Sensor
Networks (WSNs). Solutions in the Medium Access Control (MAC) layer, such as the S-MAC protocol \cite{YeEstrinToN04}, have been developed to coordinate neighboring sensors' ON-OFF schedule in order to reduce both energy consumption and packet delay. These approaches do not provide specific end-to-end latency guarantee. In \cite{WeiPaschalidis}, routing problems are considered in WSNs where each sensor switches between ON and OFF states. The authors formulate an optimization problem to pick the best path that minimizes the weighted sum of the expected energy cost and the exponent of the latency probability. In another work in \cite{NingCassSensorSleeping1}, Ning and Cassandras formulate a dynamic sleep control problem in order to reduce the energy consumed in listening to an idle channel. The idea is to sample the channel more frequently when it is likely to have traffic and less frequently when it is not. The authors extend their work in \cite{NingCassSensorSleeping2}, by formulating an optimization problem with the goal of minimizing the expected total energy consumption at the transmitter and the receiver. Dynamic programming is used to come up with an optimal policy that is shown to be more effective in cost saving than the fixed sleep time. \cite{CohenKapToN09} studies the ON-OFF scheduling in wireless mesh networks. By assuming a fixed routing tree topology used
for task transmission, each child in the tree knows exactly when its parents will wake up, and the traffic is only generated by the leaves of the tree, the authors formulate and solve an optimization problem that minimizes the total transmission energy cost while satisfying the latency and maximum energy constraints on each individual node. The major difference between this paper and the existing ones in this line of research is that we study a system with a real-time constraint for each individual task. To the best of our knowledge, ON-OFF scheduling with a real-time constraint for each individual task has not been studied extensively.

It is worth noting that there also exists papers related to the service rate control problem, in which the optimal service rate policy of either single-server or multi-server queueing systems are derived in order to minimize an average cost. A recent representative work along this line can be found in \cite{xia2017service} where Xia et al. study the service rate problem for tandem queues with power constraints. They formulate the model as a Markov decision process with constrained action space and use sensitivity-based optimization techniques to derive the conditions of optimal service rates, the optimality of the vertexes of the feasible domain for linear and concave operating cost, and an iterative algorithm that provides the optimal solution. Our problem formulation is different from these works in two aspects: \textit{(i)} We consider tasks with real-time constraints and \textit{(ii)} We include system wake-up cost on top of the service cost.

\section{System Model and Problem Formulation}
\label{Sec_model}
We consider a finite horizon scenario that a DES processes $N$ tasks with real-time constraints. In particular, task $i$, $i=1,\ldots ,N,$ has arrival time $a_{i}$ (generally random), deadline $d_{i}=a_{i}+d$, and $B$ number of operations. Both $d$ and $B$ are constants. In the \emph{off-line} setting, we assume that the task arrival time $a_{i}$ is known to the controller a priori. The DES can only operate in one of the two modes:\ ON and OFF. When it is in the OFF\ mode, there is no operating cost associated. When it is in the ON or active mode, the system can either be \emph{busy} or \emph{idling}. When the system is busy, it processes the tasks at a constant rate $R$ with fixed operating cost $C_{B}$ per unit time. When the system is idling, no tasks are waiting to be served, and the system cost is $C_{I}$ ($C_{I}$$\le$$C_{B}$) per unit time. Furthermore, we assume that whenever a transition from the OFF mode to the ON mode occurs, a fixed wake-up cost $C_{W}$ is incurred; examples of such costs include:\ the large amount of current (known as inrush current) required when a motor is turned on, the energy needed to initialize electric circuits when RF radio is turned on in a wireless device, and so on. Note that the wake-up cost may also include system wearout cost, if the system can only be turned on for certain number of times during its lifetime. In our previous work in \cite{MiaoXuWTS2015} and \cite{Miao2017ACC_ON_OFF_Control}, $C_I=C_B$. As we will show later, when $C_I$ is different from $C_B$, it does not make the analysis significantly harder, and the off-line optimal solution can still be obtained by DP. 

Our system model above is quite generic and is applicable to a wide range of engineering applications; for example, one can use ultra-low power wake-up receivers \cite{pletcher200952} to conserve energy in WSNs. Next, we formulate the off-line optimization problem. 

As we mentioned earlier, the task information is known to the controller a priori in the off-line setting. Our objective is to find the optimal ON and OFF time periods so as to \emph{(i)} finish all the tasks by their deadlines and \emph{(ii)} minimize the cost. 

\begin{definition}
	Suppose the system is woken up at $t_{1},$ put to sleep at $t_{2}$ $%
	(t_{1}<t_{2}),$ and kept active from $t_{1}$ to $t_{2}.$ Then, we call the
	time interval $[t_{1},t_{2}]$ an \textbf{Active Period} \textbf{(AP)}.
\end{definition}

\begin{definition}
	In any AP, the periods during which the system is actively serving tasks are known as \textbf{Busy Periods} \textbf{(BPs)}. The rest of the time periods in that $AP$ are known as \textbf{Idle Periods} \textbf{(IPs)}. 
\end{definition}

\begin{figure}[thpb]
	\centering
	\includegraphics[height=1.6in,width=2.8253in,angle=0]{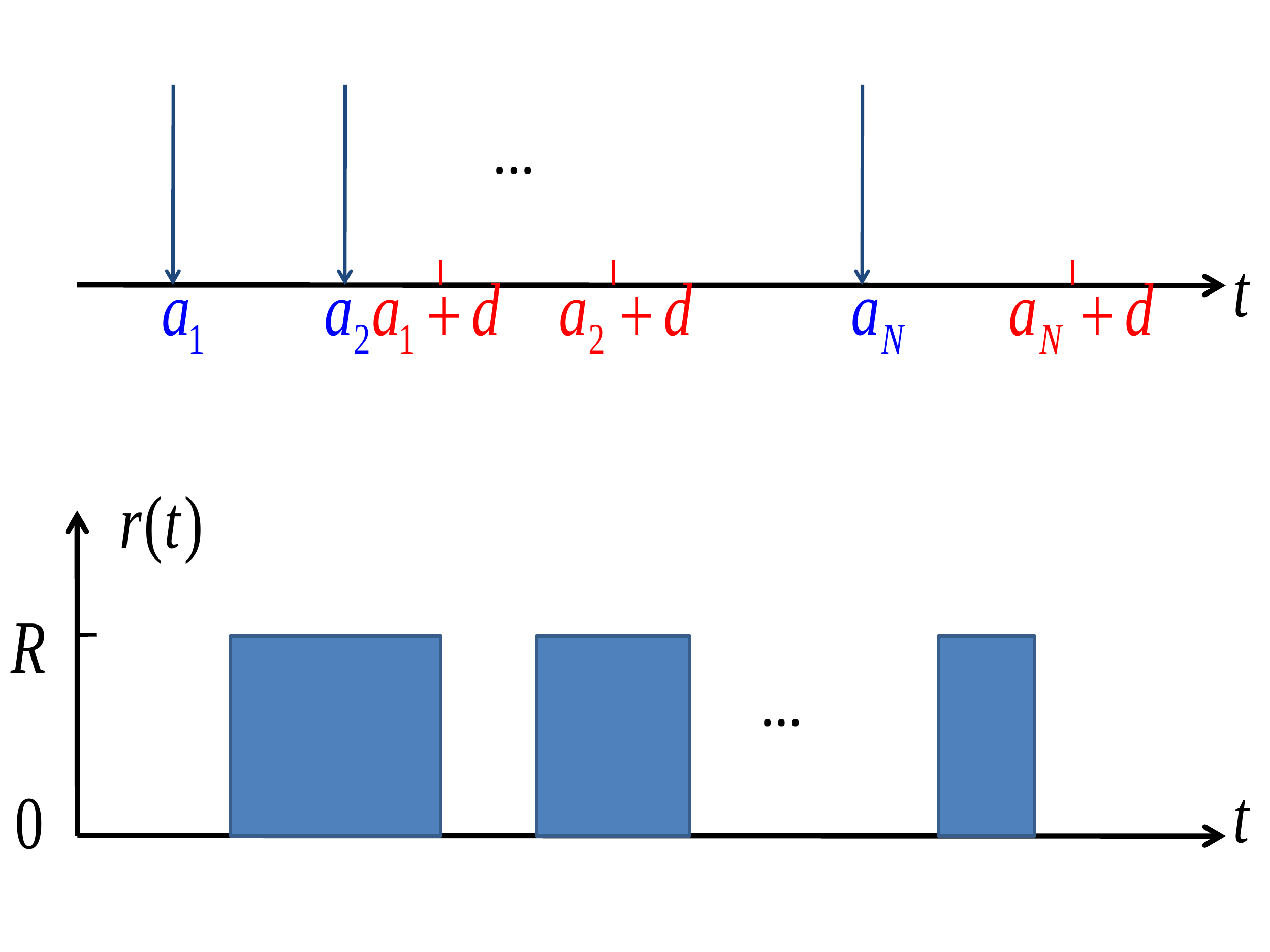}
	\caption{Off-line control illustration.}
	\label{offline_illustration}
\end{figure}

Let $r(t)$ be the rate that the system is capable of serving tasks at time $t$. It is piecewise constant and at any given time $t$, it can only be either $0$ (when the system is OFF) or $R$ (when the system is ON). See Fig. \ref{offline_illustration} for an illustration of how $r(t)$ looks like and how the APs are formed. Note that $r(t)$ is not the actual service rate since the system is only serving tasks during the \textbf{BPs}, not the \textbf{IPs}.

We now introduce the control variables. Our first control variable is $\alpha
,$ the number of APs. The second control variable is a $\alpha \times 2$
array \textbf{t} that contains $2\alpha $ time instants. These time instants
satisfy:%
\begin{equation*}
t_{i,1}<t_{i,2}<t_{j,1}<t_{j,2},\text{ }\forall i,j\in \{1,\ldots ,\alpha \},%
\text{ }i<j
\end{equation*}%
and define $\alpha $ number of APs. See Fig. \ref{offline_illustration} for illustration. The off-line problem $Q(1,N)$ can then be formulated: 
\begin{gather*}
\min_{\alpha ,\mathbf{t}}\text{ }J=\alpha C_{W}+\sum_{i=1}^{\alpha
}[C_{I}(t_{i,2}-t_{i,1}-\tau_{i,B})+C_{B}\tau_{i,B}] \\
\text{s.t. }\int_{\max (a_{j},x_{j-1})}^{x_{j}}r(t)dt=B, \\
x_{j}\leq d_{j},\text{ }x_{0}=0,\text{ }j=1,\ldots ,N \\
r(t)=R\sum_{i=1}^{\alpha }[u(t-t_{i,1})-u(t-t_{i,2})]
\end{gather*}%
where $x_{j}$ is the departure time of task $j$, $u(t)$ is the unit step function, and $\tau_{i,B}$ is the length of the busy periods in the \emph{i-th} AP. The first constraint ensures that exactly $B$ number of operations are executed for each task. The second one is the real-time constraint. The third one makes sure that the processing rate is $R$ only during each AP. Note that $\tau_{i,B}$ is dependent on the number of tasks served in $AP_{i}$. To represent $\tau_{i,B}$, we use $N_{i}^{S}$ and $N_{i}^{E}$ to denote the first (starting) task and the last (ending) task in $AP_{i}$, respectively:
\begin{gather*}
N_{i}^{E}=\underset{j\in \{1,\ldots ,N\}}{\arg \max }(d_{j}\le t_{i,2}) \\
N_{i}^{S}=\underset{j\in \{1,\ldots ,N\}}{\arg \min }(a_{j}\ge t_{i,1}) \\
\tau_{i,B}=\max((N_{i}^{E}-N_{i}^{S}+1)\frac{B}{R},0)
\end{gather*}

Notice that $Q(1,N)$ above may not always be feasible. Consider the case
that $N$ tasks arrive at the same time and need to be served in $d$
seconds. In order to meet the deadlines of all the tasks, we must have $R\geq \frac{NB}{d}$. Since $R$ is a constant, the condition above obviously is not true when $N$
is large. In this paper, we only consider the case that $Q(1,N)$\textbf{\ }%
is indeed feasible, and we have the following assumption on the task arrival
rate.

\begin{assumption}
	\label{feasibility_assumption}Within any time interval of $d$ seconds, the
	number of task arrivals must not exceed $\lfloor \frac{d}{\beta }\rfloor ,$
	where $\beta =B/R$ is the time it takes to process a single task.
\end{assumption}

We emphasize that $d$ in Assumption \ref{feasibility_assumption} is the deadline of each task upon arrival. To make the problem more interesting, we also assume that $\lfloor \frac{d}{%
	\beta }\rfloor >1.$

\begin{lemma}
	\label{Lemma_feasibility}Under Assumption \ref{feasibility_assumption}, \textbf{P1} is always feasible.
\end{lemma}

$Q(1,N)$ is a hard optimization problem, due to the nondifferentiable terms in the constraints and the objective function. It cannot be easily solved by standard optimization software. In what follows, we will first discuss optimal
off-line control, using which we will then establish the results for on-line control.

\section{Off-line Control}
\label{Sec_offline}
In this section, we focus on the off-line control problem, in
which all task arrivals are known to us a priori\textbf{.} We need to find
out when the system should wake up and start to serve the first task in
an AP. Similar to the ``just-in-time" idea exploited in \cite%
{GamNaPraUyZaInf02} for adaptive modulation, the system should wake up as
late as possible so that it may potentially reduce the active time. The
question is how late the system should wake up. This is answered by the
following results.

\begin{lemma}
	\label{Lemma_LateStartIsBetter}Suppose that tasks $\{k,\ldots ,n\}$ are all the
	tasks served in an AP on the optimal sample path of $Q(1,N)$ and starting the AP at time either $t_{A}$ or $t_{B},$ $a_{k}\leq t_{A}<t_{B}\leq
	d_{k}-\beta$, is feasible. Then, 
	\begin{equation*}
	C_{k,\ldots ,n}^{A}\geq C_{k,\ldots ,n}^{B}
	\end{equation*}
	where $C_{k,\ldots ,n}^{A}$ and $C_{k,\ldots ,n}^{B}$ are the corresponding
	costs of serving tasks $\{k,\ldots ,n\}$ in the AP for the two different starting time $t_A$ and $t_B$, respectively.
\end{lemma}

Lemma \ref{Lemma_LateStartIsBetter} indicates that an AP on the optimal
sample path of $Q(1,N)$ should be started as late as possible. We now utilize this result to figure out when
exactly the first task $k$ should be served.

\begin{lemma}
	\label{Lemma_When_to_Start_case1}If tasks $\{k,\ldots ,n\}$ are all the
	tasks served in an AP on the optimal sample path of $Q(1,N)$ and the
	number of task arrivals in $[a_{k},d_{k}-\beta )$ is $0$, then the optimal
	starting time to transmit task $k$ is $d_{k}-\beta $, i.e., 
	\begin{equation*}
	x_{k}^{\ast }=d_{k},
	\end{equation*}
where $x_{k}^{\ast }$ is the optimal departure time of task $k$.
\end{lemma}

Lemma \ref{Lemma_When_to_Start_case1} shows that we can delay the
transmission of the first task in an\ AP to $\beta $ seconds before its
deadline, provided that there are no other arrivals before that time. Next,
we discuss the case that there exists other task arrivals before $%
d_{k}-\beta $.

\begin{lemma}
	\label{Lemma_When_to_Start}Suppose task $k$ is the first task in an AP on
	the optimal sample path of $Q(1,N)$ and the number of task arrivals in $%
	[a_{k},d_{k}-\beta )$ is $m,$ $0<m\leq \lfloor \frac{d}{\beta }\rfloor -1$%
	. Let%
	\begin{equation}
	\delta _{j}=\beta (j-k)-(a_{j}-a_{k})  \label{Lemma2_1}
	\end{equation}%
	\begin{equation*}
	z=\underset{j=k+1,\ldots ,k+m}{\arg \max \{\delta _{j}}\}
	\end{equation*}%
	The optimal starting time to serve task $k$ is: 
	\begin{equation*}
	\left\{ 
	\begin{array}{cc}
	d_{k}-\beta , & \text{if }\delta _{z}\leq 0 \\ 
	d_{k}-\beta -\delta _{z}, & \text{if }\delta _{z}>0%
	\end{array}%
	\right.
	\end{equation*}
\end{lemma}

Having discussed when to wake up the system, we now find out when the system should go to sleep.
Apparently, the optimal time to end an AP depends on future task
information. In what follows, we first establish some results that identify
the end of an AP based on future task arrival information.

\begin{lemma}
	$\label{TasksApartEndAP}$If $d_{j}+C_{W}/C_{I}<a_{j+1},$ $j\in \{1,\ldots
	,N-1\},$ then task $j$ ends an \textbf{AP} on the optimal sample path of $%
	Q(1,N)$.
\end{lemma}

Lemma \ref{TasksApartEndAP} basically indicates that if the deadline of task $j$
is at least $C_{W}/C_{I}$ seconds apart from the next task arrival, then
task $j$ ends an AP on the optimal sample path. Note that this is just a
sufficient, but not necessary condition of an AP ending on the optimal
sample path. In some cases, whether a task should end an AP is determined by
not only the next arrival, but also all subsequent ones. Let $d_{0}=-\infty $
and $a_{N+1}=\infty .$ We introduce the following definition.

\begin{definition}
	Consecutive tasks $\{k,\ldots ,n\},$ $1\leq k\leq n\leq N,$ belong to a 
	\emph{super active period} (SAP) in problem $Q(1,N)$ if $%
	d_{k-1}+C_{W}/C_{I}<a_{k},$ $d_{n}+C_{W}/C_{I}<a_{n+1},$ and $%
	d_{j}+C_{W}/C_{I}\geq a_{j+1},$ $\forall j\in \{k+1,\ldots ,n-1\}.$
\end{definition}

Each SAP contains one or more APs. SAPs can be easily identified by simply
examining all the task deadlines and arrival times and applying Lemma \ref%
{TasksApartEndAP}. It implies that instead of working on the original
problem $Q(1,N)$, we now only need to focus on each SAP, which is essentially a subproblem $Q(k,n)$.

\begin{figure}[thpb]
   	\centering
   	\includegraphics[height=1.2in,width=2.12in,angle=0]{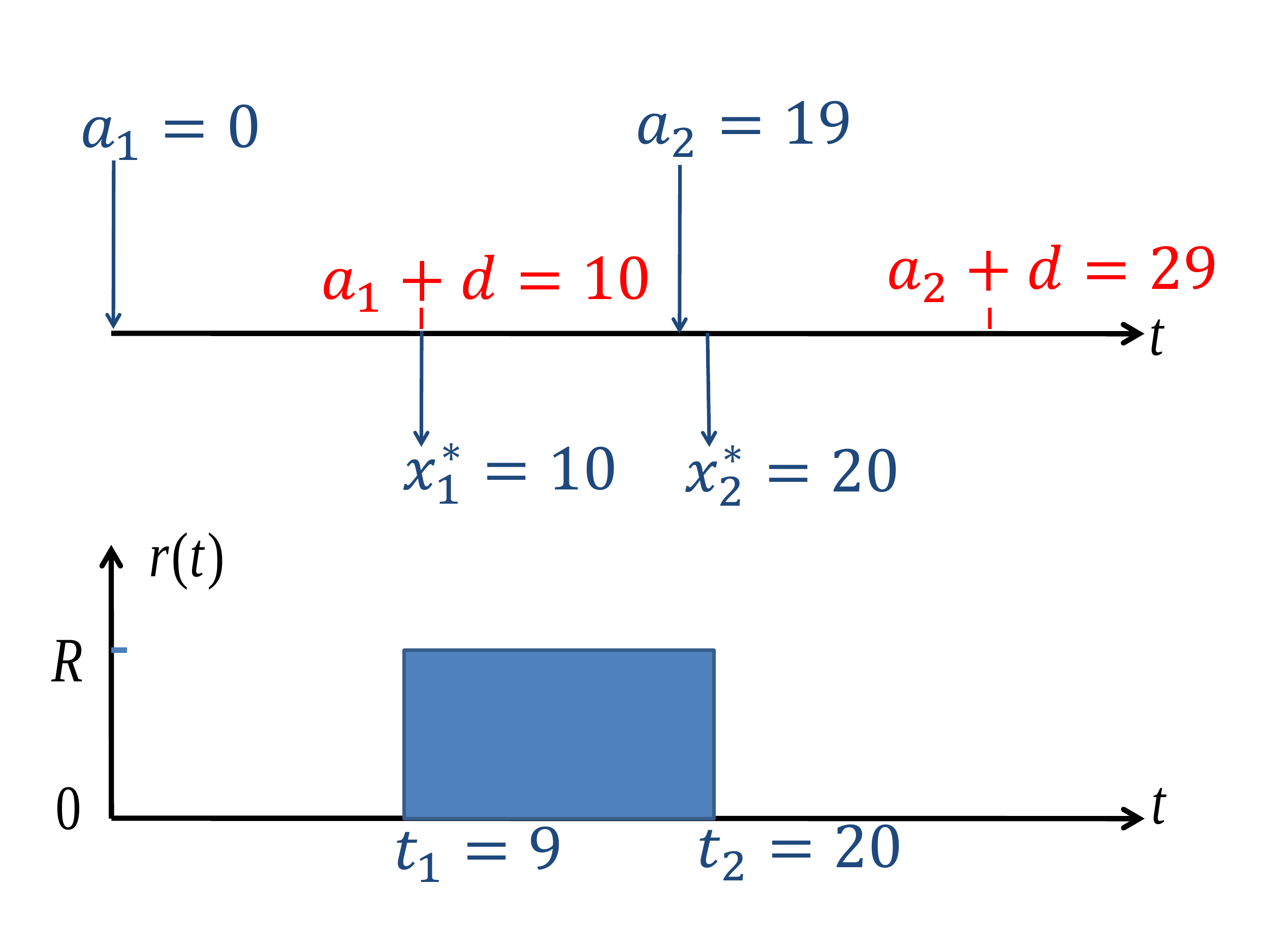}
   	\caption{Sample path \#1 of scenario \#1.}
   	\label{SP_1_Scenario_1}
   \end{figure}

\begin{figure}[thpb]
     	\centering
     	\includegraphics[height=1.2in,width=2.12in,angle=0]{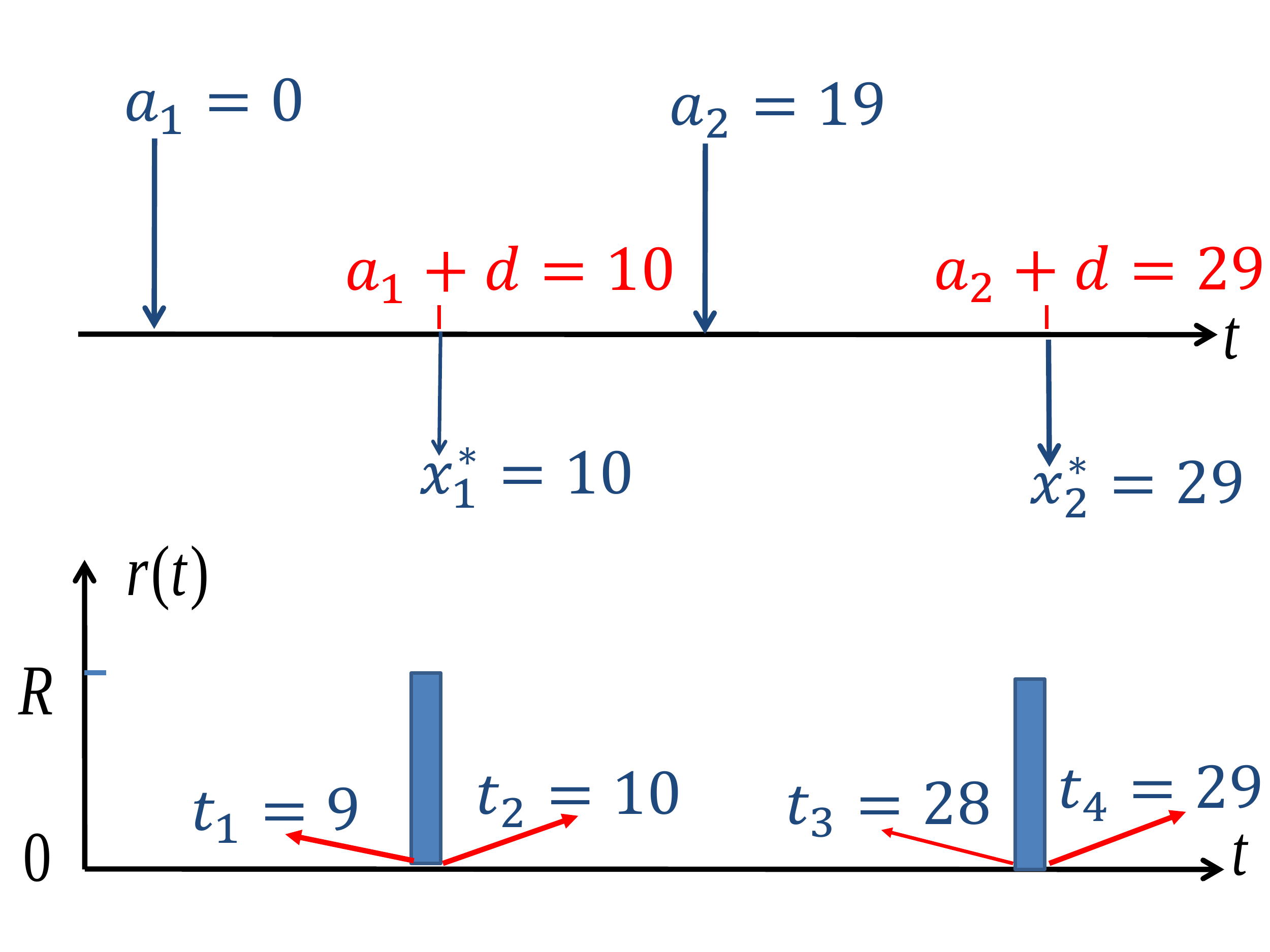}
     	\caption{Sample path \#2 of scenario \#1.}
     	\label{SP_2_Scenario_1}
  \end{figure}

We now define our decision points in each SAP. A decision point $%
x_{t}$, $t\in \{k,\ldots ,n-1\}$, is the departure time of task $t$ 
satisfies $x_{t}<a_{t+1}.$ If $x_{t} \ge a_{t+1}$, then $x_{t}$ is not a decision point because the system should stay active at $x_{t}$ and process task $t+1$. At each decision point, the control is
letting the system either go to sleep or stay awake. Let us take a look at
some examples, in which $d=10,$ $C_{W}=10,$ and $C_{B}=C_{I}=C=1.$ Note that $C_B$ and $C_I$ could be different in general; for simplicity, we let them equal to each other in the examples. We also assume
that $B=R,$ i.e., it takes a unit of time to complete a task. Fig. \ref{SP_1_Scenario_1} and Fig.  \ref{SP_2_Scenario_1}
show two different sample paths for a simple two-task scenario: $a_{1}=0$
and $a_{2}=19.$ In both sample paths, task $1$'s optimal wake up time is
determined by Lemmas 4.2 and 4.3. The only decision
point is $x_{1}$, at which the system needs to decide if it should go to
sleep or stay awake. In particular, the system in Fig. \ref{SP_1_Scenario_1} wakes up at $%
t_{1}=9 $, finishes task $1$ at its deadline $d_{1}=10,$ stays awake, and
finishes task 2 at $t_{2}=20.$ The total cost is: $%
C_{W}+C(t_{2}-t_{1})=21.$ In Fig. \ref{SP_2_Scenario_1}, the system wakes up at $t_{1}=9$,
finishes task $1$ at its deadline $d_{1}=t_{2}=10,$ and goes to sleep.
Then, it wakes up at $t_{3}=28$ (once again determined by Lemmas 4.2 and 4.3 and finishes task $2$ at $t_{4}=29.$ The total cost of this case is: $2C_{W}+C[(t_{2}-t_{1})+(t_{4}-t_{3})]=22.$ It is evident that at decision point $x_{1}=10,$ the optimal control is to let the system stay awake (shown in Fig. \ref{SP_1_Scenario_1}).

\begin{figure}[thpb]
   	\centering
   	\includegraphics[height=1.2in,width=2.12in,angle=0]{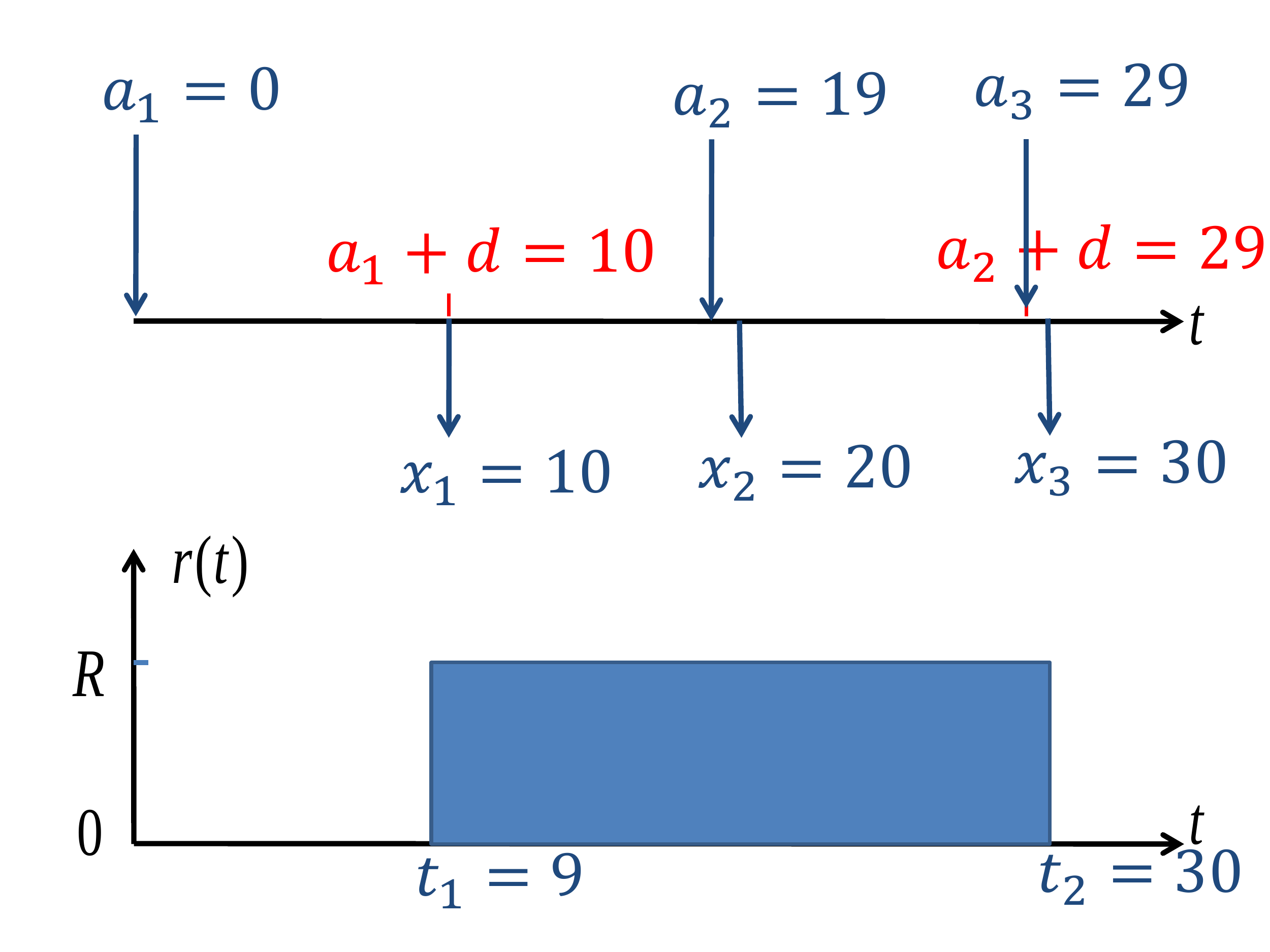}
   	\caption{Sample path \#1 of scenario \#2.}
   	\label{SP_1_Scenario_2}
   \end{figure}

\begin{figure}[thpb]
   	\centering
   	\includegraphics[height=1.2in,width=2.12in,angle=0]{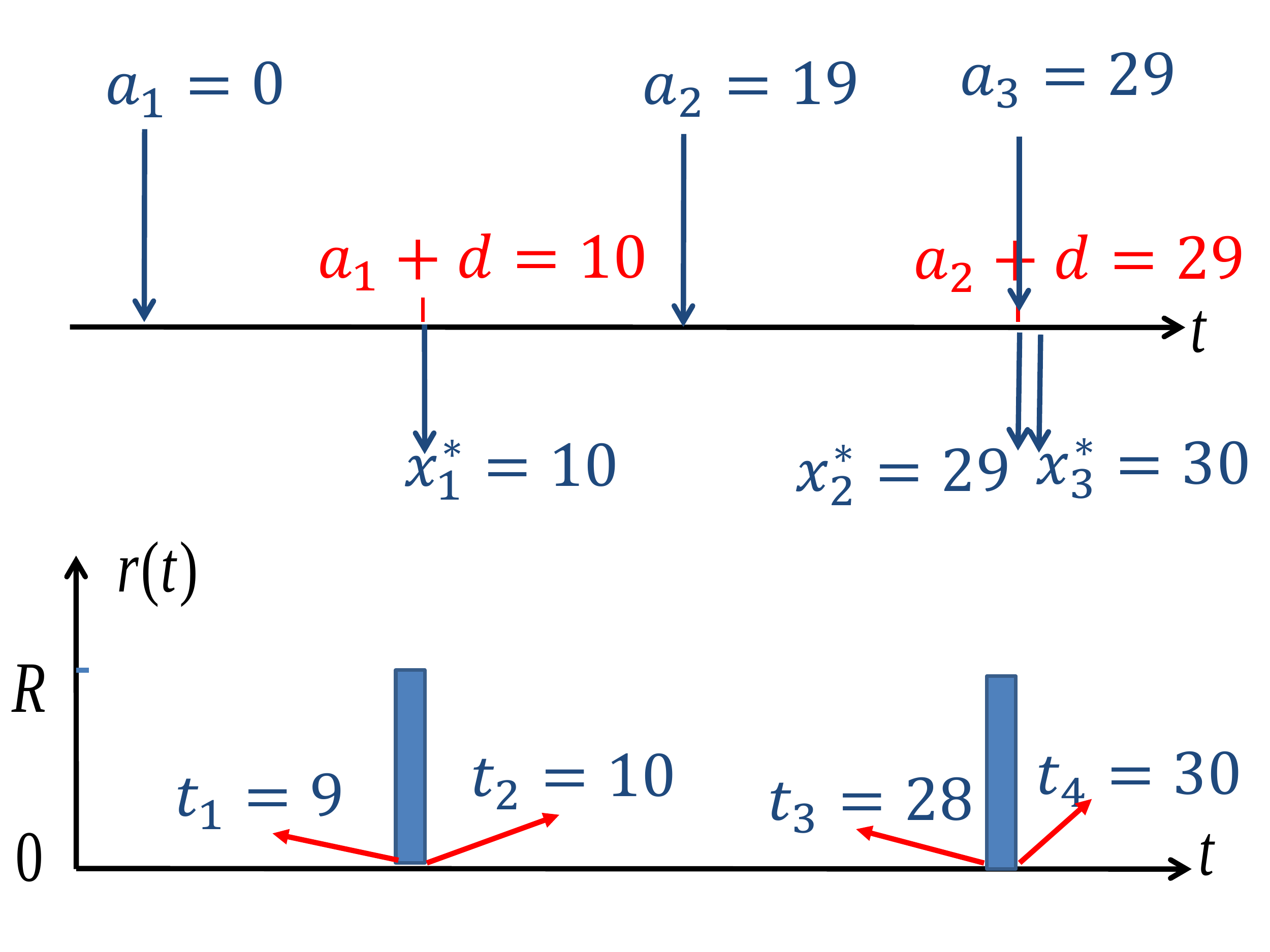}
   	\caption{Sample path \#2 of scenario \#2.}
   	\label{SP_2_Scenario_2}
   \end{figure}

Now, let us consider another scenario (Fig. \ref{SP_1_Scenario_2} and Fig. \ref{SP_2_Scenario_2}), in which we keep the previous tasks $1$ and $2$ unchanged and add task $3$. Our first decision point is again at $x_{1}=10.$ In Fig. \ref{SP_1_Scenario_2}, the system
wakes up at $t_{1}=9$, finishes task $1$ at its deadline $d_{1}=10,$ stays
awake, finishes task 2 at time $20,$ stays awake, and finally finishes task $3$ at time $t_{2}=30.$ The total cost is: $C_{W}+C(t_{2}-t_{1})=31.$ In
Fig. \ref{SP_2_Scenario_2}, the system wakes up at $t_{1}=9$, finishes task $1$ at its deadline $%
d_{1}=t_{2}=10,$ and goes to sleep. Then, it wakes up at $t_{3}=28$ and finishes tasks $2$ and $3$ at $t_{4}=30.$ The total cost of this case is: $%
2C_{W}+C[(t_{2}-t_{1})+(t_{4}-t_{3})]=23.$ It is evident that at
decision point $x_{1}=10,$ the optimal control is to let the system
go to sleep (shown in Fig. \ref{SP_2_Scenario_2}).

From the above examples, we can conclude that the optimal decision on if the
system should stay awake or go to sleep when it finishes all on-hand tasks
depends on future task arrivals (task $3$ in the examples above). A first
look at the problem seems to suggest that in the worst case, the system may have to make
a decision about if it should go to sleep or stay awake after each task
departure; the total number of possible sample paths could be as high as $2^{N}$%
, which makes the problem intractable when $N$ is large. However, a closer
look at the problem indicates that the off-line optimal ON-OFF control
problem can be solved by DP, which has been widely used to solve a large class of problems with special structural properties. In the context of DES, however, its usage has been very limited to date. For example, in \cite{MaoCTDA} and \cite{MiaoMaoCGCTONS} where the problem formulation is similar to the one in this paper, both CTDA and GCTDA algorithms are not DP-based. We will show next that for the DES studied by this paper, DP and sample path analysis can be used together to obtain the optimal solution. In particular, it is done by introducing two types of tasks: starting and following.  

\begin{definition}
	In problem $Q(k,n)$, where tasks $\{k,\ldots ,n\}$ form an SAP, the first task of any AP is called a \emph{starting}
	task. Tasks that are not starting tasks are known as \emph{following} tasks.
\end{definition}

Since the case that $k=n$ is trivial, we assume that $k<n$ in our analysis
in order to make the problem more interesting. Note that APs contain one task
only do not have following tasks. For any task $i\in\{k,\dots,n\}$, it must either be a starting
task or a following one. We are interested in finding out the optimal cost
of serving tasks $\{i,\ldots ,n\}$, and we use $Q^{S}(i,n)$ and $Q^{F}(i,n)$ to
denote the optimization problems of serving tasks $\{i,\ldots ,n\}$ when
task $i$ is a starting and following task, respective. Note that in these two problems, only tasks $\{i,\ldots,n\}$ are served and all other tasks in $\{k,\ldots,n\}$ are not considered. In problem $Q^{F}(i,n)$, the system is active when task $i$ arrives; therefore, task $i$ will be served right after its arrival.  Let $J_{i}^{S}$
and $J_{i}^{F}$ be the minimum cost of $Q^{S}(i,n)$ and $Q^{F}(i,n)$,
respectively. When $i=n$, we can easily calculate $J_{n}^{S}$ and $J_{n}^{F}:$\ $J_{n}^{S}=C_{W}+C_{B}\beta ,\text{ }J_{n}^{F}=C_{B}\beta.$ Note that $J_{n}^{F}$ does not include the wake-up cost $C_{W}$, since by assumption, task 
$n$ is a following task. The operating cost, $C_{B}\beta $, is identical in both cases. Suppose that $%
J_{i}^{S}$ and $J_{i}^{F}$ $,$ $i\in \{k+1,\ldots ,n\}$ are both known, the
next step is to find $J_{i-1}^{S}$ and $J_{i-1}^{F}.$

We first focus on $J_{i-1}^{S}.$ By assumption, task $i-1$ is a
starting task. We use Lemmas \ref{Lemma_When_to_Start_case1} and \ref{Lemma_When_to_Start} to find out the
optimal starting time of task $i-1$ in problem $Q^{S}(i-1,n).$ Let the
optimal starting time be $s_{i-1,n}^{i-1}.$ For tasks in $\{i,\ldots ,n\}$,
find task $l$ that satisfies the following:%
\begin{equation}
\label{l_for_J_S}
\begin{gathered}
s_{i-1,n}^{i-1}+(j-i+1)\beta \ge a_{j}, 
\forall j\in \{i-1,\ldots ,l-1\}, \\
\text{and }s_{i-1,n}^{i-1}+(l-i+1)\beta < a_{l}
\end{gathered}
\end{equation}
If task $l$ does not exist, then it is a trivial case that the system is
always busy serving tasks $\{i-1,\ldots ,n\}$, and there is a single AP\
that starts from $s_{i-1,n}^{i-1}$ and ends at $s_{i-1,n}^{i-1}+(n-i+2)%
\beta .$ In this case, $J_{i-1}^{S}=C_{W}+(n-i+2) \beta C_{B}$. We now consider the more interesting case that task $l$ does
exist. In particular, 
\begin{equation}
J_{i-1}^{S}=\min (V_{i-1,l}^{SS}+J_{l}^{S},V_{i-1,l}^{SF}+J_{l}^{F})
\label{J_S_minimum}
\end{equation}%
where $V_{i-1,l}^{SS}$ is the cost of serving tasks $\{i-1,\ldots ,l-1\}$
when task $l$ is a starting task:%
\begin{equation}
\label{V_SS}
V_{i-1,l}^{SS}=C_{W}+(l-i+1)\beta C_{B}
\end{equation}
$V_{i-1,l}^{SF}$ is the cost of serving tasks $\{i-1,\ldots ,l-1\}$ when task $l$
is a following task:%
\begin{equation}
\label{V_SF}
V_{i-1,l}^{SF}=C_{W}+(l-i+1)\beta C_{B}+[a_{l}-s_{i-1,n}^{i-1}-(l-i+1)\beta]C_{I}
\end{equation}

We now focus on $J_{i-1}^{F}.$ We emphasize again that in this case, task $i-1$ sees an active system upon its arrival; it will be served right away since it is the first task in $Q^{F}(i-1,n)$. For tasks in $\{i,\ldots ,n\}$, find task $l$
that satisfies the following:%
\begin{equation}
\label{l_for_J_F}
\begin{gathered}
a_{i-1}+(j-i+1)\beta \ge a_{j}, \forall j\in \{i-1,\ldots ,l-1\}, \\
\text{ and }a_{i-1}+(l-i+1)\beta < a_{l}
\end{gathered}
\end{equation}%
Once gain, task $l$ may not exist, and it corresponds to the case that the
system is always busy serving tasks $\{i-1,\ldots ,n\}.$ In this case, there
is a single AP\ that starts from $a_{i-1}$ and ends at $a_{i-1}+(n-i+2)%
\beta .$ We have $J_{i-1}^{F}=(n-i+2) \beta C_{B}$. We now consider the more interesting case that task $l$ does
exist. We have:
\begin{equation}
J_{i-1}^{F}=\min (V_{i-1,l}^{FS}+J_{l}^{S},V_{i-1,l}^{FF}+J_{l}^{F})
\label{J_F_minimum}
\end{equation}%
where $V_{i-1,l}^{FS}$ is the cost of serving tasks $\{i-1,\ldots ,l-1\}$
when task $l$ is a starting task:%
\begin{equation}
\label{V_FS}
V_{i-1,l}^{FS}=(l-i+1)\beta C_{B}
\end{equation}

$V_{i-1,l}^{FF}$ is the cost of serving tasks $\{i-1,\ldots ,l-1\}$ when task $l$
is a following task:%
\begin{equation}
\label{V_FF}
V_{i-1,l}^{FF}=(l-i+1)\beta C_{B}+[a_{l}-a_{i-1}-(l-i+1)\beta]C_{I}
\end{equation}

In Table \ref{Table_Q_k_n}, we show the algorithm that returns the optimal cost of $Q(k,n)$. This algorithm involves two more algorithms that return the optimal costs of $Q^{S}(i-1,n)$ (Table \ref{Table_Q_S}) and $Q^{F}(i-1,n)$ (Table \ref{Table_Q_F}), respectively. 

\begin{theorem}
	\label{theorem_optimal}$J_{k}^{S}$ is the optimal cost of problem $Q(k,n)$.
\end{theorem}

We have proved that when the algorithm in Table \ref{Table_Q_k_n} stops, $J_{k}^{S}$ is the optimal
cost of problem $Q(k,n)$. The corresponding optimal control, i.e., the
starting time and ending time of each AP, can be traced back iteratively by
identifying the $J_{l}^{S}$ or $J_{l}^{F}$ that each $J_{i-1}^{S}$ or $%
J_{i-1}^{F}$ points to. The procedure is provided in Table \ref{Table_Control}.

Next, we use the example in Fig. \ref{SP_1_Scenario_2} and Fig. \ref{SP_2_Scenario_2} to show how the above
algorithms work. We have three tasks $1,$ $2,$ and $3$ belong to a SAP ($k=1$
and $n=3$). Initially, $J_{n}^{S}=J_{3}^{S}=C_{B}+C\beta =11,$ and $%
J_{n}^{F}=J_{3}^{F}=C\beta =1.$ In the first iteration $(i=n=3)$, we
calculate $J_{i-1}^{S}$ and $J_{i-1}^{F}.$ To calculate $%
J_{i-1}^{S}=J_{2}^{S},$ we first figure out $s_{2,3}^{2}=28.$ Then, we find
out that no task $l$ satisfies (\ref{l_for_J_S}). Therefore, tasks $2$ and $3
$ form a single AP in problem $Q^{S}(2,3)$, and $J_{2}^{S}=12.$ To calculate 
$J_{i-1}^{F}=J_{2}^{F},$ we identify that task $l=3$ satisfies (\ref%
{l_for_J_F}). We then use (\ref{J_F_minimum}) to obtain $J_{2}^{F}=\min
(V_{2,3}^{FS}+J_{3}^{S},V_{2,3}^{FF}+J_{3}^{F})=\min
(1+J_{3}^{S},10+J_{3}^{F})=11.$ In the final iteration $(i=n-1=2)$, \ we
only need to calculate $J_{i-1}^{S}=J_{1}^{S}.$ Because $s_{1,3}^{1}=9$ and
task $l=2$ satisfies (\ref{l_for_J_S}), we use (\ref{J_S_minimum}) to
calculate $J_{1}^{S}:$ $J_{1}^{S}=\min
(V_{1,2}^{SS}+J_{2}^{S},V_{1,2}^{SF}+J_{2}^{F})=\min
(11+J_{2}^{S},20+J_{2}^{F})=23.$ This is the optimal cost obtained in Fig. \ref{SP_2_Scenario_2}. If we
follow the procedure in Table \ref{Table_Control}, we will get the exact same optimal
solution as shown in Fig. \ref{SP_2_Scenario_2}. The details are omitted.

\begin{figure}[H]
	\centering
	\includegraphics[width=0.7\textwidth]{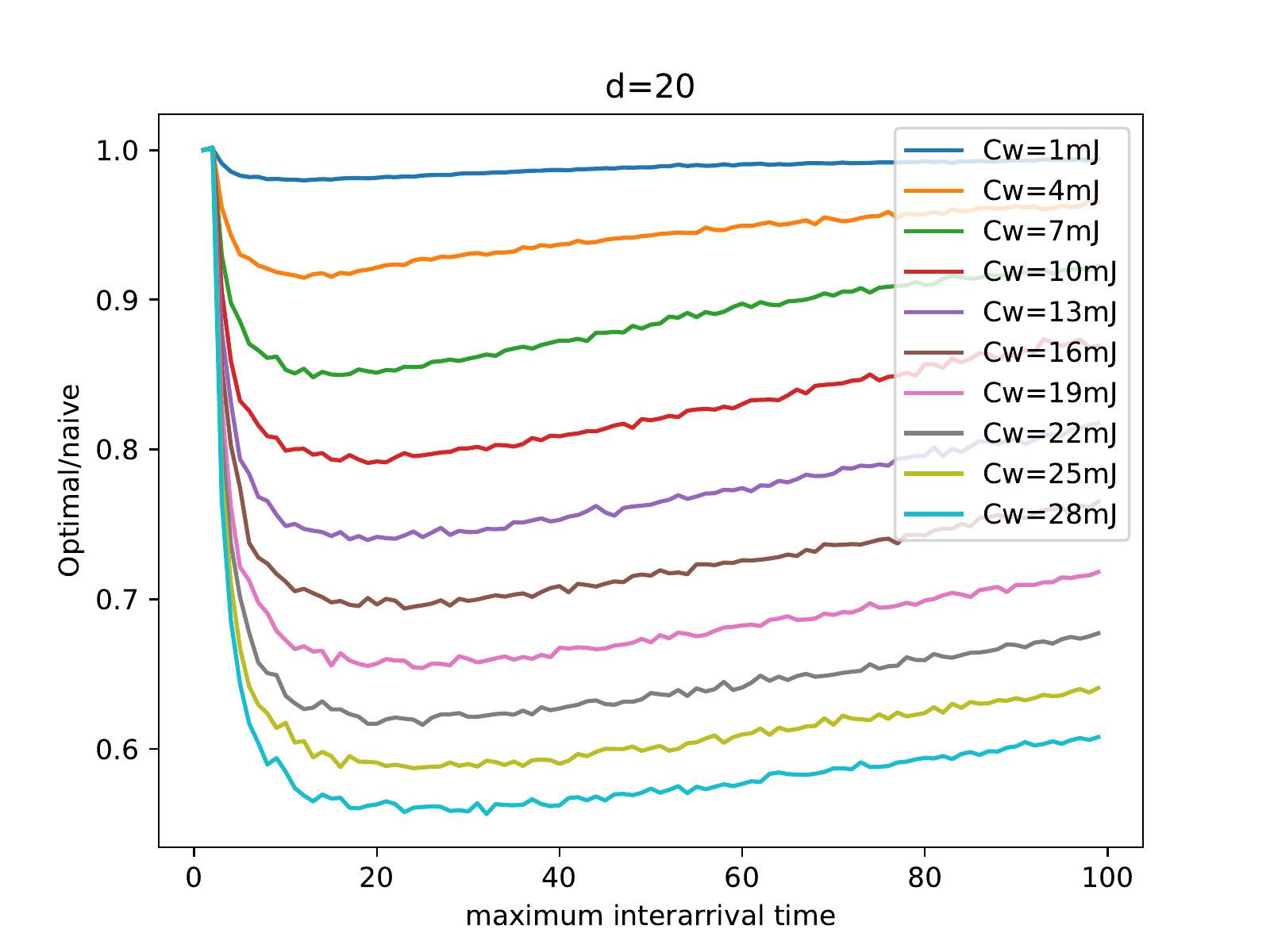}
	\caption{\label{fig:Fig4}Optimal to naive ratio under various wake-up cost and interarrival time}
	\label{Fig_optimal_naive}
\end{figure}

Next, we use simulation results to show how the optimal solution performs compared with a naive approach, in which the controller simply goes to sleep when there is no backlog and wakes up when a new task arrives. Let optimal to naive ratio be the ratio between the optimal cost and the cost of the naive controller. Fig. \ref{Fig_optimal_naive} shows how the optimal to naive ratio varies when the task arrival process and the wake-up cost $C_W$ change. In the simulation, we have $100$ runs that correspond to 100 maximum interarrival time from $1ms$ to $100ms$ with step size $1ms$. $1000$ tasks and various $C_W$ values are used in each run. The interarrival time between two adjacent tasks is uniformly distributed between 0 and the maximum interarrival time in each run. The values of the other parameters are as follows: $d=20ms$,  $C_B$=30mW, $C_I$=100$\mu$W, and $\beta=1ms$. 

We have a couple of observations. First, the cost saving of the optimal solution is greater when $C_W$ is larger. 
Second, the maximum cost saving occurs when the interarrival time is not too small or too large: when it is too small, a single AP will be sufficient to complete all the tasks, and the optimal and the naive solutions are essentially the same; when it is too large, many APs are needed, and the advantage of the optimal controller gets smaller. As we can see from the result, the cost saving of the DP algorithm in the $C_W=28mJ$ case is as large as $50\%$, and it will be ever greater when $C_W$ is higher.
\section{On-line Control}
\label{Sec_online}
In the previous section, we combine structural properties of the optimal sample path and dynamic programming to find the optimal solution to the off-line control problem. In this section, we study on-line control where future task arrival information is unknown to the controller. Essentially, the controller needs to decide the
starting time and ending time for each AP. 

\subsection{Starting an AP}

We first focus on the following questions: how can we
determine the best time to start an AP in on-line control and how different is it from the optimal time in off-line control?

\begin{figure}[thpb]
	\centering
	\includegraphics[height=2.4in,width=2.8253in,angle=0]{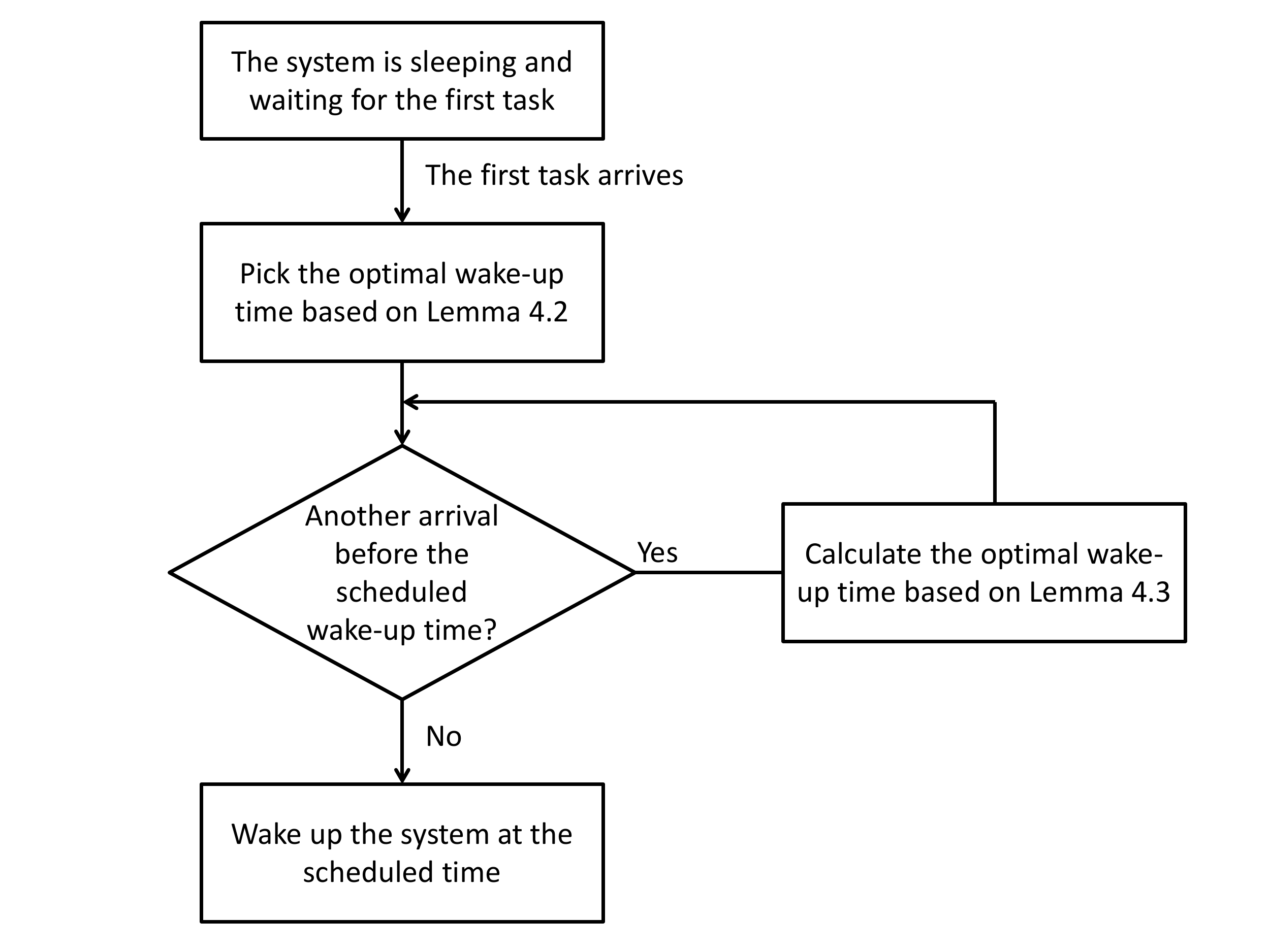}
	\caption{On-line control: starting an AP}
	\label{start_of_AP}
\end{figure}

Fig. \ref{start_of_AP} shows the proposed on-line control mechanism for
determining the wake-up time. It is an iterative algorithm that dynamically
adjusts the wake-up time based upon the backlog and the newly available task
information. Initially, right after the first task arrives, the scheduled
wake-up time is $a_{1}+d-\beta$ (determined by Lemma \ref{Lemma_When_to_Start_case1}). If there are other task arrivals before
the scheduled wake-up time, the controller will recalculate the wake-up time
using the results in Lemma \ref{Lemma_When_to_Start}; otherwise, the system
will be woken up at the scheduled time.

\begin{lemma}
	\label{decision_after_arrival}Suppose that tasks $\{k,\ldots ,n\}$ form an AP on the optimal sample path of $Q(1,N)$. If the
	system is OFF before task $k$ arrives in on-line control, then
	the wake-up time returned by the on-line control mechanism in Fig. \ref%
	{start_of_AP} is optimal.
\end{lemma}

Lemma \ref{decision_after_arrival} indicates that for on-line control, the lack
of future task information does not incur any penalty when starting an AP:\
the optimal time to start an AP can be determined iteratively using the
backlog and the newly available task information.

We now turn our attention to ending an AP in on-line control. 

\subsection{Ending an AP}
When all backlogged tasks have been served in an on-line setting, the controller needs to decide when to end an AP and put the system to sleep. This decision depends on future task information and the values of idling cost $C_I$ and wake-up cost $C_W$. For example, if the next task $t+1$ arrives very soon, the optimal control at decision point $x_{t}$ might be letting the system stay active; conversely, if the next task $t+1$ arrives after a long time, then the system perhaps should go to sleep at decision point $x_{t}$. 

When some future task information is known, techniques such as Receding Horizon Control (RHC) can be utilized to make decisions. In this paper, we focus on the scenario that no future task information is available at all.

In general, the control at each decision point is the following: let the system stay awake for another $\theta_{t}$ seconds. If no task arrives within the $\theta_{t}$ seconds, then put the system to sleep after $\theta_{t}$ seconds; o.w., serve the newly arrived tasks and wait for the next decision point. Note that the subscript $t$ indicates $\theta_{t}$ could be different at each decision point. 

Let $J^{*}$ be the optimal cost of the off-line problem $Q(1,N)$ and $\widetilde{J}$ be the cost of the on-line controller. Our objective is to develop competitive on-line controllers which can quantify their worst-case performance deviation from the optimal off-line solution. 

One challenge of competitive analysis is to find out the worst-case scenario. In our problem, the unnecessary cost in on-line control occurs when the system is idling: the controller must decide if and when to sleep. Therefore,     the worst case occurs when each AP contains only one task so that the decision has to be made over and over again for every single task. This property actually simplifies our analysis, and in particular, we tackle the competitive ratio problem from two different aspects: a deterministic controller and a randomized one.   
\subsubsection{Deterministic Controller}
We first consider a deterministic controller in which $\theta_{t}$ is a fixed constant value $\theta$. The on-line controller is \emph{c-competitive} if $\widetilde{J}(I,\theta)\le c J^{*}, \forall I\in\mathscr{I}$,
where $\mathscr{I}$ is the set of all possible task arrival instances and $I$ is one task arrival instance. $c$ is called the competitive ratio of the deterministic on-line controller and is essentially the \emph{upper-bound (i.e., worst case)} of the ratio between the on-line cost $\widetilde{J}$ and the off-line optimal cost $J^{*}$.

\begin{lemma}
\label{Lemma_deterministic_competitive_ratio}
The best competitive ratio $c^*$ of a deterministic controller is obtained when $\theta=C_W/C_I$, and $\lim_{N\to\infty} c^*=(2+\gamma)/(1+\gamma)$, where $N$ is the number of tasks and $\gamma=C_B\beta/C_W$.
\end{lemma}

Lemma \ref{Lemma_deterministic_competitive_ratio} shows that the competition ratio of a deterministic algorithm depends on the ratio between $C_B\beta$, the cost of serving one task, and $C_W$, the cost of waking up the system. If this ratio is very small, then the competitive ratio is close to $2$; if the ratio is very large, then the competitive ratio is close to $1$.

\subsubsection{Randomized Controller} In a different methodology, we assume that $\theta_{t}$ is determined by a randomized algorithm that returns a value based on certain probability distribution $P$. During on-line control, the controller essentially is playing a game with an adversary (i.e., the task arrival process). Our job is to find the optimal probability distribution and the corresponding competitive ratio. We point out that the competitive ratio of a randomized on-line algorithm $A$ is defined with respect to a specific type of adversary. In this paper, we assume an oblivious adversary \cite{ben1994power}, in which the worst instance for the randomized algorithm $A$ is chosen without the the knowledge of the realization of the random variable used by $A$. We say randomized algorithm $A$ is \emph{c-competitive} if $E_{P}[\widetilde{J}(A,I)] \le cJ^{*}(I), \forall I \in \mathscr{I}$, 
where $\widetilde{J}(A,I)$ is the cost of algorithm $A$ under task arrival instance $I$ in on-line control and $J^{*}(I)$ is the corresponding off-line optimal cost.
Note that the task arrival instance $I$ must be fixed before the expectation is taken.
The competition ratio of randomized algorithm $A_{P}$ (algorithm $A$ using probability distribution $P$) is:
\begin{equation*}
c(A_{P})=\underset{I\in\mathscr{I}}{\sup}\frac{E_{P}[\widetilde{J}(A_{P},I)]}{J^{*}(I)}
\end{equation*} 
Our goal is to find the best possible probability distribution that yields the best competitive ratio $c^{*}$: 
\begin{equation*}
\label{minimax}
c^{*}=\underset{P}{\inf}\text{ }\underset{I\in\mathscr{I}}{\sup}\frac{E_{P}[\widetilde{J}(A_{P},I)]}{J^{*}(I)}
\end{equation*}
This is essentially a minimax problem, and one way of solving it is to use Yao's minimax principle \cite{yao1977probabilistic}, which states: a randomized algorithm may be viewed as a random choice between deterministic algorithms; in particular, the competitive ratio of a randomized algorithm against any oblivious adversary is the same as that of the best deterministic algorithm under the worst-case distribution of the adversary's input. In our case, the adversary's input is the task arrival instance after each AP. Let its probability distribution be $G$. Using Yao's principle and von Neumann minimax theorem, we get:
\begin{equation}
\label{maxmin}
c^{*}=\underset{G}{\sup}\text{ }\underset{A\in\mathscr{A}}{\inf}\frac{E_{G}[\widetilde{J}(A,I_{G})]}{J^{*}(I_{G})}
\end{equation}
where $\mathscr{A}$ is the set of all randomized algorithms, $I_{G}$ is a specific task arrival instance under probability distribution $G$, and the expectation is now performed with respect to $G$. We now use the following lemma to find $c^*$.

\begin{lemma}
	\label{lemma_competitive_random}
	The best competitive ratio $c^*$ of a randomized controller is obtained when $\theta_t$ is a random variable $X$, whose probability density function is \begin{eqnarray}
	\label{f_x}
	f_X(x) &=&\left\{ 
	\begin{array}{cc}
	\frac{1}{\frac{C_W}{C_I}(e-1)}e^{x/(C_W/C_I)}, & \text{if }x\leq C_W/C_I \\ 
	0, & \text{if }x>C_W/C_I%
	\end{array}%
	\right. 
	\end{eqnarray}
When this controller is used, $\lim_{N\to\infty} c^*=(\gamma+1.58)/(\gamma+1)$, where $\gamma=C_B\beta/C_W$.
\end{lemma}

Lemma \ref{lemma_competitive_random} shows that the competition ratio of a random controller also depends on the ratio between $C_B\beta$ and $C_W$. If this ratio is very small, then the competitive ratio is close to $1.58$; if the ratio is very large, then the competitive ratio is close to $1$.

\section{Conclusions}
\label{Sec_conclusions}
In this paper, we study the optimal ON-OFF control problem for a class of DESs with real-time
constraints. The DESs have operating costs $C_{B}$ and $C_{I}$ per unit time and wake-up cost $C_{W}$. Our goal is to switch the system between the ON and the OFF states so as to minimize cost and satisfy real-time constraints. In particular, we consider a homogeneous case that all tasks
have the same number of operations and each one's deadline is $d$ seconds after
the arrival time. For the off-line scenario that all task information is known to the controller
a priori, we show that the optimal solution can be obtained via a two-fold
decomposition: $(i)$ super active periods that contain one or more active
periods can be identified easily using the task arrival times and deadlines
and $(ii)$ the optimal solution to each super active period can be solved
using dynamic programming. Simulation results show that compared with a simple heuristic, the cost saving of the DP algorithm can be 50\% or more.

In on-line control, we show that the best time to start an AP can be obtained via an iterative algorithm and is guaranteed to be the same as the off-line problem. To decide the best time to end an AP in the on-line setting where no future task arrival information is available, we evaluate both deterministic and  random controllers and derive their competitive ratios; these results quantify the worst-case on-line performance deviation from the optimal off-line solution. 

%\begin{acknowledgements}
%\end{acknowledgements}



\appendix 

\begin{center}
	APPENDIX
\end{center}

\textbf{Proof of Lemma \ref{Lemma_feasibility}:} To see this, consider the solution that the system is woken
up at $a_{1}$ and stays active until $d_{N}.$ Because in any $d$ seconds, 
\begin{equation*}
Rd\geq \lfloor \frac{d}{\beta }\rfloor B,
\end{equation*}%
that is, the number of task departures is not less than the number of task
arrivals, the backlog is always zero at the integer multiples of $d$
seconds. This means that all task arrivals can be served within $d$
seconds. Hence, the proposed solution is always feasible under Assumption %
\ref{feasibility_assumption}. $\blacksquare $

\textbf{Proof of Lemma \ref{Lemma_LateStartIsBetter}:}\ Because%
\begin{equation*}
a_{k}\leq t_{A}<t_{B}\leq d_{k}-\beta
\end{equation*}%
we have 
\begin{equation}
x_{j}^{A}\leq x_{j}^{B},\text{ }j=k,\ldots ,n  \label{departure_early}
\end{equation}%
where $x_{k}^{A},\ldots ,x_{n}^{A}$ and $x_{k}^{B},\ldots ,x_{n}^{B}$ are
the task departure times in the two sample paths, respectively. It means
that the departure time of task $j$ in sample path $A$ is no later than that
in sample path $B$. Note that (\ref{departure_early}) holds because the
sensor stays on in the AP and $t_{A}<t_{B}.$

Let $t_{EA}$ and $t_{EB}$ be the ending time of the AP when the starting time is $t_A$ and $t_B$, respectively. We have:
\begin{equation*}
C_{k,\ldots
	,n}^{A}=C_{w}+C_{B}(n-k+1)\beta+C_I[t_{EA}-t_A-(n-k+1)\beta]
\end{equation*}
\begin{equation*}
C_{k,\ldots ,n}^{B}=C_{w}+C_{B}(n-k+1)\beta+C_I[t_{EB}-t_B-(n-k+1)\beta]
\end{equation*}
where the three terms in each equation correspond to the wake-up cost, the cost of serving the $n-k+1$ tasks, and the idling cost, respectively.

By assumption that $t_A<t_B$, we know for sure that the idling time of sample path $A$ is not less than that of sample path $B$, i.e., $t_{EA}-t_A \ge t_{EB}-t_B$. Therefore, 
\begin{equation*}
C_{k,\ldots ,n}^{A}\geq C_{k,\ldots ,n}^{B} \text{  }\blacksquare
\end{equation*}

\textbf{Proof of Lemma \ref{Lemma_When_to_Start_case1}:}\ Time $d_{k}-\beta $ is the latest time to serve task $k$%
, in order to meet its hard deadline requirement. Invoking Lemma \ref%
{Lemma_LateStartIsBetter}, we only need to show that $Q(1,N)$ is still
feasible when we delay the transmission of task $k$ until $d_{k}-\beta .$
Under Assumption \ref{feasibility_assumption} and in the worst case, there
could be $\lfloor \frac{d}{\beta }\rfloor -1$ tasks $\{k+1,\ldots
,k+\lfloor \frac{d}{\beta }\rfloor -1\}$ arriving at $d_{k}-\beta .$ It
means that at time $d_{k}-\beta ,$ we have $\lfloor \frac{d}{\beta }%
\rfloor $ tasks in the backlog. If we start transmitting all these tasks $%
\{k,\ldots ,k+\lfloor \frac{d}{\beta }\rfloor -1\}$ at $d_{k}-\beta $,
then it takes maximum $d$ seconds to send all of them, and each task's
deadline is met. Again, under Assumption \ref{feasibility_assumption}, the
earliest time that task $k+\lfloor \frac{d}{\beta }\rfloor $ can arrive is $%
d_{k}.$ If the system stays active at $d_{k}+d-\beta ,$ then task $%
a_{k+\lfloor \frac{d}{\beta }\rfloor }$ can be served by its deadline $%
d_{k}+d.$ Similarly, all subsequent tasks can be transmitted by their
deadlines. Therefore, $Q(1,N)$ is still feasible after we postpone task $%
k $'s transmission time to $d_{k}-\beta .$ $\blacksquare $

\textbf{Proof of Lemma \ref{Lemma_When_to_Start}:} Invoking Lemma \ref{Lemma_LateStartIsBetter} again, we need
to show that $d_{k}-\beta $ and $d_{k}-\beta -\delta $ are the latest
feasible starting times for the two cases, respectively.

In the worse case, there could be $\lfloor \frac{d}{\beta }\rfloor -m-1$
tasks arriving at $d_{k}-\beta .$ As we have shown in Lemma \ref%
{Lemma_When_to_Start_case1}, these tasks and all subsequent ones can be
served before their deadlines as long as we start the AP no later than $%
d_{k}-\beta .$ Therefore, we only need to focus on the tasks that arrive
before $d_{k}-\beta .$

\emph{Case 1:}$\ \delta _{z}\leq 0.$ This implies that 
\begin{equation}
\frac{a_{j}-a_{k}}{j-k}\geq \beta ,\text{ for }j\in \{k+1,\ldots ,k+m\}
\label{All_greater_than_beta}
\end{equation}

$d_{k}-\beta $ is the latest possible starting time for task $k.$ We need
to show that $Q(1,N)$ is still feasible when we start serving task $%
k $ at $d_{k}-\beta ,$ i.e., starting serving tasks at this time will
satisfy the real-time constraints for tasks $\{k+1,\ldots ,k+m\}$. If task $%
k $ is done at $d_{k}-\beta ,$ then the departure time $x_{j}$ of
task $j$, $j\in \{k+1,\ldots ,k+m\},$ is 
\begin{equation*}
x_{j}=d_{k}+\beta (j-k)
\end{equation*}%
From (\ref{All_greater_than_beta}), we have%
\begin{equation*}
x_{j}\leq d_{k}+a_{j}-a_{k}=a_{j}+d=d_{j}
\end{equation*}%
Thus, the deadlines of all the tasks $\{k+1,\ldots ,k+m\}$ are met, and $%
d_{k}-\beta $ is the optimal starting time.

\emph{Case 2:} $\delta _{z}>0$. We need to show $d_{k}-\beta -\delta _{z}$
is a feasible starting time for all tasks $\{k,\ldots ,k+m\}.$ We first show
causality. The starting time $s_{j}$ of task $j\in
\{k+1,\ldots ,k+m\}$ is:$\ $ 
\begin{equation}
s_{j}=d_{k}-\beta -\delta _{z}+\beta (j-k)  \label{Lemma2_3}
\end{equation}

Using (\ref{Lemma2_1}),%
\begin{eqnarray}
s_{j} &=&d_{k}-\beta -[\beta (z-k)-(a_{z}-a_{k})]+\beta (j-k)  \label{Sj}
\\
&=&d_{k}-\beta +\beta (j-z)+(a_{z}-a_{k})  \notag \\
&=&d-\beta +\beta (j-z)+a_{z}  \notag
\end{eqnarray}

We use contradiction to prove it. Suppose $s_{j}<a_{j},$ we have 
\begin{equation*}
d-\beta +\beta (j-z)+a_{z}<a_{j}\text{, i.e.,}
\end{equation*}%
\begin{equation}
d-\beta +\beta (j-z)<a_{j}-a_{z}  \label{inequality}
\end{equation}

1) When $k<j\leq z\leq k+m,$ we have 
\begin{equation}
a_{j}-a_{z}\leq 0  \label{Inequality_1}
\end{equation}%
By Assumption \ref{feasibility_assumption}, 
\begin{equation*}
\beta (j-z)=-\beta (z-j)\geq -(d-\beta ),
\end{equation*}%
i.e., 
\begin{equation}
d-\beta +\beta (j-z)\geq 0  \label{Inequality_2}
\end{equation}%
Combining (\ref{Inequality_1}) and (\ref{Inequality_2}), (\ref{inequality})
is not true.

2) When $k<z<j\leq k+m$, we have 
\begin{equation*}
d-\beta >a_{j}-a_{z}>0
\end{equation*}%
and 
\begin{equation*}
d-\beta +\beta (j-z)>d
\end{equation*}%
Combining the two inequalities, we conclude that (\ref{inequality}) is not
true either.

We can now assert 
\begin{equation*}
s_{j}\geq a_{j}
\end{equation*}%
which satisfies causality.

Next, we show the departure time of each task $j\in \{k+1,\ldots ,k+m\}$ is
before the task's deadline. Again, we use $x_{j}$ to denote the departure
time of task $j,$ and 
\begin{equation*}
x_{j}=s_{j}+\beta
\end{equation*}%
Invoking (\ref{Sj}), 
\begin{equation*}
x_{j}=d+\beta (j-z)+a_{z}
\end{equation*}%
We need to show 
\begin{equation*}
x_{j}=d+\beta (j-z)+a_{z}\leq a_{j}+d,
\end{equation*}%
i.e., 
\begin{equation}
\beta (j-z)+a_{z}\leq a_{j}  \label{proof_feasibility}
\end{equation}

From (\ref{Lemma2_1}), we have 
\begin{eqnarray*}
	\delta _{j} &=&\beta (j-k)-(a_{j}-a_{k})\leq \\
	\delta _{z} &=&\beta (z-k)-(a_{z}-a_{k}),\text{ } \\
	j &=&k+1,...,k+m
\end{eqnarray*}%
Rearranging the terms above, we obtain (\ref{proof_feasibility}).

Finally, the departure time of task $z$ is exactly $a_{z}+d,$ indicating
that $d_{k}-\beta -\delta _{z}\ $is the latest possible time to start
serving task $z$. $\blacksquare $

\textbf{Proof of Lemma \ref{TasksApartEndAP}:} We use $x_{j}^{\ast }$ and $s_{j+1}^{\ast }$ to denote the
departure time of task $j$ and the starting time of task $j+1,$
respectively, on the optimal sample path of $Q(1,N)$. Using Lemma \ref{Lemma_feasibility}, we have 
\begin{equation}
x_{j}^{\ast }\leq d_{j}  \label{x<d}
\end{equation}%
From casualty, 
\begin{equation}
s_{j+1}^{\ast }\geq a_{j+1}  \label{s>a}
\end{equation}%
By assumption, we have 
\begin{equation}
a_{j+1}-d_{j}>C_{W}/C_{I}  \label{a-d>cw/ca}
\end{equation}%
Combining (\ref{x<d}), (\ref{s>a}), and (\ref{a-d>cw/ca}), we get%
\begin{equation}
s_{j+1}^{\ast }-x_{j}^{\ast }>C_{W}/C_{I}  \label{s-x>cw/ca}
\end{equation}%
Next, we use a contradiction argument to prove the lemma. Let the optimal
sample path of $Q(1,N)$ be $sp^{\ast }$ and the corresponding cost is $%
J^{\ast }$. Suppose that task $j$ does not end an \textbf{AP} on $sp^{\ast }$%
. It means that the system stays active from $x_{j}^{\ast }$ to $%
s_{j+1}^{\ast }.$ The optimal cost is then $J^{\ast }=(s_{j+1}^{\ast }-x_{j}^{\ast })C_{I}+J_{R},$ 
where $J_{R}$ is the rest of the cost beyond time interval $[x_{j}^{\ast
},s_{j+1}^{\ast }].$ Consider another sample path $sp^{^{\prime }}$, which
is identical to $sp^{\ast }$, except that the system goes to sleep at $%
x_{j}^{\ast }$ and wakes up at $s_{j+1}^{\ast }.$ The system cost is now $J^{^{\prime }}=C_{W}+J_{R}.$
Using (\ref{s-x>cw/ca}), we obtain $J^{^{\prime }}<J^{\ast },$ which
contradicts the assumption that $sp^{\ast }$ is the optimal sample path. $%
\blacksquare $

\textbf{Proof of Theorem \ref{theorem_optimal}:} We use induction to prove it.

\emph{Step 1}: Task $n$ can either be a starting task or a following task.
When it is a starting task, it is obvious that $J_{n}^{S}$ is the optimal
cost of $Q^{S}(n,n).$ When it is a following task, it is also obvious that $%
J_{n}^{F}$ is the optimal cost of $Q^{F}(n,n).$

\emph{Step 2}: Suppose that $J_{j}^{S}$ is the optimal cost of problem $%
Q_{j}^{S}(j,n),$ and $J_{j}^{F}$ is the optimal cost of problem $%
Q_{j}^{F}(j,n),$ $j\in \{i,\ldots ,n\}.$ We need to show that $J_{i-1}^{S}$
and $J_{i-1}^{F}$ are the optimal cost of problem $Q_{i-1}^{S}(i-1,n)$ and $%
Q_{i-1}^{F}(i-1,n),$ respectively$.$ Since the proofs are similar, we only
show that $J_{i-1}^{S}$ is the optimal cost of problem $Q_{i-1}^{S}(i-1,n).$
By assumption, task $i-1$ is a starting task. We can use Lemmas \ref{Lemma_When_to_Start_case1} and \ref{Lemma_When_to_Start} to find $s_{i-1,n}^{i-1},$ the optimal
starting time of task $i-1$. We now discuss two cases:

Case 1:\ Task $l$ that satisfies (\ref{l_for_J_S}) does not exist.

It implies that 
$s_{i-1,n}^{i-1}+(j-i+1)\beta >a_{j},\forall j\in \{i-1,\ldots ,n\},$ i.e., the system is busy serving tasks whenever a task $j\in \{i-1,\ldots
,n\}$ arrives. Therefore, there is no reason to go to sleep, and tasks $%
\{i-1,\ldots ,n\}$ form a single AP. From Line 14 of Table \ref{Table_Q_S}, $J_{i-1}^{S}=C_{W}+(n-i+2) \beta C_{B}$ is the optimal cost of problem $Q^{S}(i-1,n)$.

Case 2:\ Task $l$ that satisfies (\ref{l_for_J_S}) does exist.

In this case, task $l$ has not arrived when task $l-1$ departs the system.
It has two subcases:\ the system should either go to sleep when task $l-1$
departs or stay awake (and serve task $l$ when it arrives). The subcase that
yields a smaller cost is the optimal solution, and this is calculated in (%
\ref{J_S_minimum}). $\blacksquare $

\textbf{Proof of Lemma \ref{decision_after_arrival}:} We consider two cases.

\emph{Case 1:}\ The optimal wake-up time is $d_{k}-\beta .$ This happens
when either Lemma \ref{Lemma_When_to_Start_case1} or the $\delta _{z}\leq 0$
case of Lemma \ref{Lemma_When_to_Start} applies. The on-line control
mechanism picks the same wake-up time upon the arrival of task $k$, and it
does not change. Therefore, the wake-up time in on-line control is the same
as the optimal wake-up time on the optimal sample path.

\emph{Case 2:} The optimal wake-up time is $d_{k}-\beta -\delta _{z}.$ This
happens when the $\delta _{z}<0$ case of Lemma \ref{Lemma_When_to_Start}
applies. In on-line control, the initial wake-up time is set to $d_{k}-\beta
.$ With the arrival of tasks between $a_{k}$ and $d_{k}-\beta ,$ this
scheduled time is adjusted to $d_{k}-\beta -\delta _{j}$, for some $j\in
\{k+1,\ldots ,k+m\}.$ By definition of $\delta _{z},$ we have 
\begin{eqnarray*}
	d_{k}-\beta -\delta _{j} &\geq &d_{k}-\beta -\delta _{z} \\
	&=&d_{k}-\beta -[\beta (z-k)-(a_{z}-a_{k})]\text{ } \\
	&=&d-\beta -\beta (z-k)+a_{z} \\
	&\geq &a_{z}
\end{eqnarray*}%
The above implies that all intermediate wake-up times and the optimal
wake-up time are after the arrival of task $z$. Therefore, the on-line
control policy is able to wake up the system at the optimal time $%
d_{k}-\beta -\delta _{z}$ after task $z$ arrives. $\blacksquare $

\textbf{Proof of Lemma \ref{Lemma_deterministic_competitive_ratio}:} The worst-case happens when each AP only contains a single task. After each task is served, the system stays active for $\theta$ seconds and then goes to sleep; it wakes up again after the next task arrives. For any $\theta$, we have the ratio between the on-line cost and the optimal cost:
\begin{equation}
\label{c_of_theta}
c(\theta)=\frac{C_W+NC_B\beta+(N-1)(C_I\theta+C_W)}{C_W+NC_B\beta+(N-1)\min(C_I\theta,C_W)}
\end{equation} 
where the numerator is the on-line cost and the denominator is the off-line cost. Note that both costs have three terms: the first term $C_W$ is the wake-up cost for serving the very first task; the second term $NC_B\beta$ is the actual cost of serving the $N$ tasks; and the last term is the cost between two adjacent tasks.
We can rewrite (\ref{c_of_theta}) into:
\begin{equation*}
c(\theta)=\frac{C_W/N+C_B\beta+(N-1)(C_I\theta+C_W)/N}{C_W/N+C_B\beta+(N-1)\min(C_I\theta,C_W)/N}
\end{equation*}
It follows that
\begin{equation*}
\lim_{N\to\infty}c(\theta)=\frac{C_B\beta+(C_I\theta+C_W)}{C_B\beta+\min(C_I\theta,C_W)}.
\end{equation*}
Because 
\begin{equation*}
\frac{C_I\theta+C_W}{\min(C_I\theta,C_W)} \ge 2 
\end{equation*}
and the equality holds when $\theta=C_W/C_I$, we have
\begin{equation*}
c^*=\lim_{N\to\infty}c(C_W/C_I)=\frac{C_B\beta+2C_W}{C_B\beta+C_W}=\frac{2+\gamma}{1+\gamma}\text{ } \blacksquare
\end{equation*}

\textbf{Proof of Lemma \ref{lemma_competitive_random}:} Similar to the deterministic algorithm case, the worst case also occurs when each AP only contains a single task. At the $i$-th decision point, the system stays active for $\theta_t=X$ seconds, where $X$ is a random variable returned by algorithm $A$, and then goes to sleep if no tasks arrive during this period. For serving $N$ tasks, the ratio between the on-line cost and the optimal-cost is:
\begin{equation}
\label{c_of_theta_random}
c(\theta_t)=\frac{C_W+NC_B\beta+(N-1)E_G[\widetilde{J}_b(A,I_G)]}{C_W+NC_B\beta+(N-1)J_b^*(I_G)}
\end{equation} 
where the numerator is the on-line cost and the denominator is the off-line cost. Note that both costs have three terms: the first term $C_W$ is the wake-up cost for serving the very first task; the second term $NC_B\beta$ is the actual cost of serving the $N$ tasks; and the last term is the cost between two adjacent tasks. Note that the reason why the expectation is taken with respect to $G$ is due to the insight provided by equation (\ref{maxmin}).
We can rewrite (\ref{c_of_theta_random}) into:
\begin{equation*}
c(\theta_t)=\frac{C_W/N+C_B\beta+(N-1)E_G[\widetilde{J}_b(A,I_G)]/N}{C_W/N+C_B\beta+(N-1)J_b^*(I_G)/N}.
\end{equation*}
It follows that
\begin{equation}
\label{lim_random}
\lim_{N\to\infty}c(\theta_t)=\frac{C_B\beta+E_G[\widetilde{J}_b(A,I_G)]}{C_B\beta+J_b^*(I_G)}.
\end{equation}

Let $y$ be the time it takes for the next task to arrive after the system finishes serving the previous task. Similar to other on-line scheduling scenarios such as the ski rental and the snoopy caching problems \cite{karlin1994competitive} , it can be seen via variational analysis that $E_{G}[\widetilde{J}_b(A,I_{G})]/J_b^{*}(I_{G})$ is uniform with respect to $y$, i.e., it is independent from $y$. Letting $\tilde{c}=E_{G}[\widetilde{J}_b(A,I_{G})]/J_b^{*}(I_{G})$, our goal is to come up the best algorithm $A$ that minimizes $\tilde{c}$. It has been shown in \cite{karlin1994competitive} that $\tilde{c}^*$ is $e/(e-1)\approx 1.58$, and the probability distribution $P$ that achieves this ratio is in (\ref{f_x}). Note that  $J_b^*(I_G)=\min(C_I\times y,C_W)$. The impact of $\tilde{c}$ to (\ref{lim_random}) is the greatest when $J_b^*(I_G)$ takes the maximum value $C_W$. In this case, $E_G[\widetilde{J}_b(A,I_G)]$ takes its value $1.58C_W$. Therefore, 
\begin{equation*}
\lim_{N\to\infty} c^*=\frac{C_B\beta+1.58C_W}{C_B\beta+C_W}=\frac{\gamma+1.58}{\gamma+1}\text{ } \blacksquare
\end{equation*}

\begin{table}[h]
	\centering%
	\begin{tabular}{ll}
		1. & $J_{n}^{S}=C_{W}+C_{B}\beta ,\text{ }J_{n}^{F}=C_{B}\beta ,\text{ and 
		}$ \\ 
		& $\text{set both }J_{n}^{S}\rightarrow next\text{ and }J_{n}^{F}\rightarrow
		next\text{ to NULL.}$ \\ 
		2. & $for$ $(i=n;i-k>=1;i--)$ $\{$ \\ 
		3. & \ \ Initialize $J_{i-1}^{S}\rightarrow next\text{ and }%
		J_{i-1}^{F}\rightarrow next$ to NULL \\ 
		4. & \ \ Solve $Q^{S}(i-1,n)$ \\ 
		5. & \ \ Solve $Q^{F}(i-1,n)$ \\ 
		6. & \}%
	\end{tabular}
	\caption{The algorithm that returns the optimal cost of $Q(k,n)$}\label{Table_Q_k_n}
\end{table}

\begin{table}[h]
	\centering%
	\begin{tabular}{ll}
		1. & Use Lemmas \ref{Lemma_When_to_Start_case1} and \ref{Lemma_When_to_Start} to find \\ 
		& $s_{i-1,n}^{i-1},$ the optimal starting time of task $i-1$ \\ 
		2. & If (there exists $l$ that satisfies (\ref{l_for_J_S})) \{ \\ 
		3. & \ \ \ Calculate $V_{i-1,l}^{SS}$ and $V_{i-1,l}^{SF}$ using (\ref{V_SS}) and (\ref{V_SF}), respectively \\ 
		4. & \ \ \ If ($V_{i-1,l}^{SS}+J_{l}^{S}\leq V_{i-1,l}^{SF}+J_{l}^{F})$ \{
		\\ 
		5. & $\ \ \ \ \ \ J_{i-1}^{S}=V_{i-1,l}^{SS}+J_{l}^{S}$ \\ 
		6. & \ \ \ \ \ \ $J_{i-1}^{S}\rightarrow next$ $=$ $J_{l}^{S}$ \\ 
		7. & \ \ \ \} \\ 
		8. & \ \ \ else \{ \\ 
		9. & $\ \ \ \ \ \ J_{i-1}^{S}=V_{i-1,l}^{SF}+J_{l}^{F}$ \\ 
		10. & \ \ \ \ \ \ $J_{i-1}^{S}\rightarrow next=$ $J_{l}^{F}$ \\ 
		11. & \ \ \ \} \\ 
		12. & \} \\ 
		13. & else \{ // single AP case \\ 
		14. & \ \ \ $J_{i-1}^{S}=C_{W}+(n-i+2) \beta C_{B}$ \\ 
		15. & \}%
	\end{tabular}
	\caption{The algorithm that returns the optimal cost of $Q^{S}(i-1,n)$}\label{Table_Q_S}
\end{table}

\begin{table}[h]
	\centering%
	\begin{tabular}{ll}
		1. & If (there exists\ task $l$ that satisfies (\ref{l_for_J_F})) \{ \\ 
		2. & \ Calculate $V_{i-1,l}^{FS}$ and $V_{i-1,l}^{FF}$ using (\ref{V_FS}) and (\ref{V_FF}), respectively\\ 
		3. & \ \ If ($V_{i-1,l}^{FS}+J_{l}^{S}\leq V_{i-1,l}^{FF}+J_{l}^{F})$ \{ \\ 
		4. & $\ \ \ \ J_{i-1}^{F}=V_{i-1,l}^{FS}+J_{l}^{S}$ \\ 
		5. & \ \ \ \ $J_{i-1}^{F}\rightarrow next=$ $J_{l}^{S}$ \\ 
		6. & \ \} \\ 
		7. & \ else \{ \\ 
		8. & $\ \ \ \ J_{i-1}^{F}=V_{i-1,l}^{FF}+J_{l}^{F}$ \\ 
		9. & \ \ \ $J_{i-1}^{F}\rightarrow next= $ $J_{l}^{F}$ \\ 
		10. & \ \} \\ 
		11. & \} \\ 
		12. & else \{ //single AP case \\ 
		13. & \ \ \ $J_{i=1}^{F}=(n-i+2) \beta C_{B}$ \\ 
		14. & \}%
	\end{tabular}
	\caption{The algorithm that returns the optimal cost of $Q^{F}(i-1,n)$}\label{Table_Q_F}
\end{table}

\begin{table}[h]
	\centering
	\begin{tabular}{ll}
		1. & $J=J_{k}^{S},$ $i=J.task=J^{\prime }s$ subscript, and \\ 
		& $J.type=J^{\prime }s$ superscript \\ 
		2. & while ( $J\rightarrow next$ is not NULL)\{ \\ 
		3. & \ \ $J^{\prime }=J->next$ \\ 
		4. & \ $\ next\_task=J^{\prime }.task$ \\ 
		5. & \ \ $next\_type=J^{\prime }.type$ \\ 
		6. & \ \ If $(J.type=``S")$\{ \\ 
		7. & \ \ \ \ \ AP starts at $s_{i,n}^{i};$ \\ 
		8. & \ \ \} \\ 
		9. & \ \ If ($next\_type=``S"$) \{ \\ 
		10. & \ \ \ \ \ AP ends after task $next\_task-1$ is served; \\ 
		11. & \ \ \ \ \ $J=J^{\prime }$ and $i=J.task$; continue; \\ 
		12. & \ \ \} \\ 
		13. & \ \ If ($next\_type=``F"$) \{ \\ 
		14. & \ \ \ \ \ Keep the system active through $a_{next\_task}$ \\ 
		15. & \ \ \} \\ 
		16. & \ \ $J=J^{\prime }$ and $i=J.task$ \\ 
		17. & \}%
	\end{tabular}
	\caption{The procedure that returns the optimal control to $Q(k,n)$}\label{Table_Control}
\end{table}

\end{document}